\documentclass{amsart}

\usepackage{epsfig}
\usepackage{amsmath}
\usepackage{amsfonts}
\usepackage{amssymb}
\usepackage{amsthm}

\def\R{\mathbb{R}}
\def\x{{\bf x}}
\def\y{{\bf y}}

\def\varphi0{{\cal Q}}

\newcommand{\be}{\begin{equation}}
\newcommand{\ee}{\end{equation}}
\newcommand{\baa}{\begin{array}}
\newcommand{\eaa}{\end{array}}
\newcommand{\ba}{\begin{eqnarray}}
\newcommand{\ea}{\end{eqnarray}}

\newtheorem{thm}{\textsc{Theorem}}[section]
\newtheorem{lem}[thm]{\textsc{Lemma}}
\newtheorem{pro}[thm]{\textsc{Proposition}}

\newtheorem{definition}[thm]{\textsc{Definition}}

\title[]{The Singular Set of Higher Dimensional Unstable Obstacle Type Problems}
\author[J. Andersson]{John Andersson}
\address{Mathematics Institute,
University of Warwick
Coventry CV4 7AL, UK}
\email{j.e.andersson@warwick.ac.uk}
\author[H. Shahgholian ]{Henrik Shahgholian}
\address{Department of Mathematics, Royal Institute of Technology,
100~44  Stockholm, Sweden}
\email{henriksh@math.kth.se}
\urladdr{http://www.math.kth.se/~henriksh/}
\author[Georg Weiss]{Georg S. Weiss}
\address{Department of Mathematics, Heinrich Heine University, 40225 D\"usseldorf}
\email{weiss@math.uni-duesseldorf.de}
\thanks{$2000$ {\it Mathematics Subject Classification.\/} Primary
35R35, Secondary 35B40, 35J60.}
\thanks{H. Shahgholian has been supported in part by
the Swedish Research Council.
Both J. Andersson and G.S. Weiss thank the G\"oran Gustafsson Foundation
for visiting appointments to
KTH}
\date{}
\begin{document}
    \maketitle

\tableofcontents

\section{Introduction.}

In this paper we will investigate the singular points of the following
unstable free boundary problem:
\begin{equation}\label{Eq}
\Delta u= -\chi_{\{ u>0\}} \quad\quad\textrm{ in } B_1(0)
\end{equation}
where $\chi_{\{u>0\}}$ is the characteristic function of the set $\{u>0\}$.

This problem was first investigated by G.S Weiss and R. Monneau \cite{MW}.
In \cite{MW}, $C^{1,1}$-regularity locally energy minimising and maximal 
solutions of (\ref{Eq}) is shown. There is also some 
discussion regarding the possibility of the existence of singular
points, that is points $x^0\in B_1(0)$ such that $u\notin C^{1,1}(B_r(x^0))$
for any $r>0$. Such points are proved to be totally unstable \cite{MW}.

Let us formally define {\sl singular points} before we proceed.

\begin{definition}\label{Sdef}
Let $u$ be a solution to (\ref{Eq}). Then we define $S(u)$, the set of singular
points of $u$, according to
$$
S(u)=\left\{\x\in B_1(0);\; u\notin C^{1,1}(B_r(\x)) \textrm{ for any }r>0 \right\}.
$$
Furthermore we will denote by $S_{n-2}(u)$ the singular points of co-dimension 2:
$$
S_{n-2}(u)=\bigg\{\y\in S(u);\; \lim_{r_j\to 0}\frac{u(r_j\x +\y)}{\|u(r_j \x+\y)\|_{L^2(B_1(0))}}= 
$$
$$
=Q\circ \frac{x_{n-1}^2-x_n^2}{\|x_{n-1}^2-x_n^2\|_{L^2(B_1(0))}} \textrm{ for some }Q\in \mathcal{Q} \textrm{ and }r_j\to 0\bigg\}
$$
where $\mathcal{Q}$ is the matrix group of rotations of $\R^n$.
\end{definition}

It was shown in \cite{MW} or \cite{ASW3} that if $y\in S(u)$ then
$$
\lim_{r_j\to 0}\frac{u(r_j\x +\y)}{\|u(r_j\x+\y)\|_{L^2(B_1(0))}}\in \mathcal{P}_2
$$
if the right hand side is defined, here 
$\mathcal{P}_2$ is the set of homogeneous second order harmonic
polynomials of degree 2. Since the only homogeneous second order harmonic 
polynomial, up to translations, rotations and multiplicative constants, 
in $\R^2$ is $x_1^2-x_2^2$ it follows that $S_{n-2}$ singles out the 
singular points with co-dimension 2 singularities.

In \cite{AW} two of the authors showed rigorously that 
singular points exists, that is there exist a solution $u$ to (\ref{Eq}) such
that $S(u)\ne\emptyset$. This investigation was followed by the
authors in \cite{ASW2} and \cite{ASW3} where we investigated and provided 
a total classification of singular points in $\R^2$ and $\R^3$ respectively.

In this paper we 
intend to prove that in $\R^n$ the singular points of smallest co-dimension
are locally contained in a $C^1-$manifold of dimension $n-2$ and
that the free boundary $\Gamma_u$, defined
$$
\Gamma_u=\{\x\in B_1(0);\; u(\x)=0\},
$$
consists of two $C^1$ manifolds of dimension
$n-1$ intersecting orthogonally at such singular points.

Our main theorem is
\begin{thm}\label{Main}
Let $u$ be a solution to (\ref{Eq}) and assume that
\begin{equation}\label{Pining}
\lim_{r_j\to 0}\frac{u(r_j\x)}{\|u(r_j \x)\|_{L^2(B_1)}}=\frac{x_{n-1}^2-x_n^2}{\| x_{n-1}^2-x_n^2  \|_{L^2(B_1)}}
\end{equation}
for some sequence $r_j\to 0$ (In particular, $0\in S_{n-2}(u)$). Then 
\begin{equation}\label{1stclaim}
\lim_{r\to 0}\frac{u(r_j\x)}{\|u(r_j \x)\|_{L^2(B_1)}}=\frac{x_{n-1}^2-x_n^2}{\| x_{n-1}^2-x_n^2  \|_{L^2(B_1)}}
\end{equation}
and for each $\eta>0$ there exists an $r_\eta>0$ such that
$$
S\cap B_{r_0}(0)\cap \left\{ \x;\; \sum_{i=1}^{n-2}x_i^2\le \eta(x_{n - 1}^2+x_n^2)  \right\}
$$
consists of two $C^1$ manifolds intersecting at right angles at the origin.

Furthermore there is a constant $r_0(u)>0$ such that the set 
$$
S_{n-2}=\bigg\{ \y;\; u(\y)=|\nabla u(\y)|=0 \textrm{ and } \lim_{r\to 0}\frac{u(r\x +\y)}{\|u(r \x+\y)\|_{L^2(B_1(0))}}= 
$$
$$
=Q\circ \frac{x_{n-1}^2-x_n^2}{\|x_{n-1}^2-x_n^2\|_{L^2(B_1(0))}} \textrm{ for some }q\in \mathcal{Q} \bigg\}
$$
is contained in a $C^1$ manifold of dimension $(n-2)$ in $B_{r_0}(0)$. 
\end{thm}

We would like to place this result in a long tradition of 
regularity result for parametric non-linear PDE. In particular
we may view the free boundary $\Gamma_u=\{\x\in B_1(0);\; u(\x)=0\}$
as a parametric surface with singular points in $S(u)$.

Some of the most famous result in this area are the results by Bombieri, De 
Giorgi, Giusti and Simmons (\cite{BGG}, \cite{simons}) that states that
no minimal cones exists for minimal surfaces in $n<8$. We should also 
mention the result by B. White \cite{White} where uniqueness of tangent cones
for 2-dimensional minimal surfaces is proved. From our point of view White's 
proof is interesting in that he uses a Fourier series expansion in constructing
comparison surfaces. However, we work in $n-$dimensions which means that
our Fourier expansions are considerably more subtle and involved than
those that appear in \cite{White}.

Singularities in parametric problems have appeared in other areas of mathematics
as well and our results have some similarities to the theory for
harmonic mappings (\cite{SimonBook} for a good overview). One could
also mention a certain similarity with the theory of singularities that
arise for $\alpha$-uniform measures \cite{KP}.

Equation (\ref{Eq}) also arises in several applications
for instance in solid combustion (see the references in \cite{MW}),
the composite membrane problem (\cite{chanillo1},
\cite{chanillo2}, \cite{blank}, \cite{shahgholian}, \cite{chanillokenig}, \cite{chanillokenigto}), climatology (\cite{diaz})
and fluid dynamics (\cite{struwe}).

Our proof will be based on a dynamic systems approach where we 
project a solution $\frac{u(r\x)}{r^2}$ onto the harmonic second order
polynomials, call this projection $\Pi(u,r,0)$ (see Definition 
\ref{projection}). By a careful analysis of the PDE we will be able to
estimate $\Pi(u,r,0)-\Pi(u,r/2,0)$. Close to a singular point
we have that $\frac{u(r\x)}{r^2}\approx \Pi(u,r,0)+Z_{\Pi(u,r,0)}$ where
$$
\begin{array}{ll}
\Delta Z_{\Pi(u,r,0)}=-\chi_{\{\Pi(u,r,0)>0\}} & \textrm{ in }\R^n \\
Z_{\Pi(u,r,0)}(0)=|\nabla Z_{\Pi(u,r,0)}(0)|=0 & \\
\lim_{|\x|\to \infty}\frac{Z_{\Pi(u,r,0)}(\x)}{|\x|^3}=0 & \\
\Pi(Z_{\Pi(u,r,0)},1,0)=0. & 
\end{array}
$$

If we disregard lower order terms we may consider the map
$\mathcal{F}(\Pi(u,r,0))= \Pi(u,r/2,0)$ defined by
$$
\mathcal{F}(\Pi(u,r,0))=\Pi(u,r,0)+\Pi(Z_{\Pi(u,r,0)},1/2,0).
$$
The blow-up is unique
if $\lim_{k\to \infty}\mathcal{F}^k(\Pi(u,r,0))$ exists.

Since the harmonic second order polynomials form a finite dimensional
space. The map $\mathcal{F}$ is a map between finite dimensional vector spaces.
The main difficulty is that $\mathcal{F}$ is highly non-linear and
we need quite subtle estimates to characterise the map. On the
positive side we may write down $\Pi(u,r,0)$ explicitly, modulo lower order
terms,
by means of Theorem \ref{KarpM} by Karp and Margulis \cite{KarpMargulis}.
The definition of $\mathcal{F}$ involves a Fourier series expansion
of $-\chi_{\Pi(u,r,0)}$ on the unit sphere. Our main effort will be to 
estimate the Fourier coefficients in this expansion when $\Pi(u,r,0)/ \sup_{B_1}|\Pi| \approx
x_{n-1}^2-x_n^2$. For further
details on the idea of the proof we refer the reader to \cite{ASW3}.

\section{List of notation:}
\begin{enumerate}
\item $\delta$ will denote a vector in $\R^{n-2}$, we will always assume that $|\delta|< < 1$. We also define $\tilde{\delta}=\sum_{i=1}^{n-2}\delta_i$.

\item $\mathcal{P}_2$ will denote the second order homogeneous polynomials.

\item $S(u)$ and $S_{n-2}(u)$ are the singular set and the singular set 
of co-dimension 2 respectively, defined in Definition \ref{Sdef}.

\item The mapping $F$ is defined in equation (\ref{Fdef}).

\item $\Pi(u,r,x^0)$ is defined in Definition \ref{projection}.

\item The average of $u$ in $\Omega$ will be denoted $(u)_{\Omega}$.

\item By $dA$ we mean an area element of the surface under considration.

\item We will use Landau's $O(r)$ notation to indicate a term that is bounded 
from by $Cr$ for a universal constant $C$. That is $f(x)=O(r)$ if and
only if $|f(x)|\le Cr$ for a universal constant $C$. Similarly,
$f(r)\ge O(r)$ means that $f(r)\ge Cr$ for some universal constant $C>0$
etc.

\item $p_\delta(x)=\sum_{i=1}^{n-2}\delta_i x_i^2+(1-\tilde{\delta})x_{n-1}^2-x_n^2$, in particular $p_0(x)=x_{n-1}^2-x_n^2$.

\item $Z_{p_\delta}$ is defined in (\ref{Zdef}).

\item $\mathcal{Q}$ is the matrix-group of rotations of $\R^n$.

\item The functions $B_i(\delta)$, $B(\delta)$, $C_i(\delta)$ 
and $C(\delta)$ are defined in (\ref{Bi}), (\ref{B}), Proposition \ref{Cjest}
and the remark after that Proposition respectively.

\end{enumerate}

\section{Background Material and General Strategy.}

In this section we will state some of the results of \cite{ASW3}
and outline our strategy (which is similar to the strategy of \cite{ASW3}).

Our starting observation is the following proposition (Proposition 5.1 in 
\cite{MW})

\begin{pro}\label{fixedcenter}
Let $u$ be a solution of (\ref{Eq}) in $B_1(0)$
and let us consider a point $\x^0\in S(u).$
Then 
$$
\lim_{r_j\to 0}\frac{u(r_j\x +\x^0)}{\|u(r_j x+\x^0)\|_{L^2(B_1(0))}}\in \mathcal{P}_2
$$
for each sequence $r_j\to 0$ such that the limit exists.
\end{pro}
The proof is a fairly standard application of a monotonicity formula.

If $u$ is a solution to (\ref{Eq}) then $\Delta u\in L^\infty$ which 
directly implies that $D^2 u\in BMO(B_{1/2}(0))$ which in particular 
implies, via the Sobolev inequality, that for $\x^0\in S(u)\cap B_{1/2}(0)$
\begin{equation}\label{BMO}
\frac{u(r\x +\x^0)-\frac{1}{2} ( \x-\x^0) (D^2 u)_{B_r(\x^0)} (\x-\x^0)}{r^2}
\end{equation}
is locally bounded in $L^2$ and  pre-compact. It will be convenient 
for some calculations later to subtract a harmonic polynomial in (\ref{BMO})
instead of 
the polynomial $\frac{1}{2}( \x-\x^0) (D^2 u)_{B_r(\x^0)} (\x-\x^0)$. We make the following definition.

\begin{definition}\label{projection}
By $\Pi(u,r,\x^0)$ we will denote the projection operator onto $\mathcal{P}_2$ defined as follows:
$\Pi(u,r,\x^0)=\tau_r p$, where $\tau_r\in \mathbb{R}^+$ and $p\in \mathcal{P}_2$ satisfies
$\sup_{B_1}|p|=1$ as well as
$$
\inf_{h\in  \mathcal{P}_2}\int_{B_1(0)}\Big|D^2\left(\frac{u(r\x+\x^0)}{r^2}\right)-D^2 h \Big|^2=
\int_{B_1(0)}\Big|D^2\left(\frac{u(r\x+\x^0)}{r^2}\right)-\tau D^2 p \Big|^2.
$$

We will often write $\Pi(u,r)$ when $\x^0$ is either the origin or given by
the context. By definition $\tau_r=\sup_{B_1}|\Pi(u,r)|$
and $p_r=\Pi(u,r)/\tau_r$.
\end{definition}

It is a simple consequence of the $BMO$ estimate (\ref{BMO}) that if 
$\x^0\in S(u)\cap B_{1/2}(u)$ then (Proposition 3.7 in \cite{ASW3})
\begin{equation}\label{c1alphaBMO}
\left\| \frac{u(r\x +\x^0)}{r^2}-\Pi(u,r,\x^0)\right\|_{C^{1,\alpha}(B_1)}\le C_\alpha(\sup_{B_1}|u|,n).
\end{equation}

If $\x^0\in S(u)$ then 
\begin{equation}\label{loggrowth}
\sup_{B_r(\x^0)}|u|> cr^2\ln(1/r)
\end{equation}
for $0<r<r_0(u,\x^0)$ and some small $c>0$. To be more 
precise it is known that (c.f. Lemma 5.1 in \cite{ASW3}). 

\begin{lem}\label{superquad} Let $u$ be a solution to (\ref{Eq}) in $B_1$
such that $\sup_{B_1}|u|\le M$ and
$u(0)=|\nabla u(0)|=0$.
Then there exist $\rho_0>0$ and $r_0>0$ such that if 
\begin{equation}\label{delta0}
\sup_{B_1}|\Pi(u,r)|\ge \frac{1}{\rho_0}
\end{equation}
for an $r\le r_0$ then
$$
\sup_{B_1}|\Pi(u,r/2)|>\sup_{B_1}|\Pi(u,r)|+\eta_0/2,
$$
where $\eta_0$ is a universal constant.
\end{lem}
The Lemma is proved 
for $n=3$ in \cite{ASW3} but the proof is the same in arbitrary dimension.

This estimate together with (\ref{c1alphaBMO}) implies that
$u(\cdot+\x^0)=\Pi(u,r,\x^0)+$a lower order perturbation. 
Using the pre-compactness in $C^{1,\alpha}$ (c.f. Equation (\ref{c1alphaBMO})) of 
\begin{equation}
 \frac{u(r_j\x +\x^0)}{r_j^2}-\Pi(u,r_j,\x^0)
\end{equation}
for some sequence $r_j\to 0$ we may extract a sub-sequence, which we still
denote by $r_j$, such that
$$
\lim_{r_j\to 0}\left(\frac{u(r_j\x +\x^0)}{r_j^2}-\Pi(u,r_j,\x^0) \right)
=Z_p(\x)
$$ 
for some function $Z_p$. It is not difficult to see that
$Z_p$ is the unique solution to
\begin{equation}
\begin{array}{ll}\label{Zdef}
\Delta Z_p=-\chi_{\{p(\x)>0\}} & \textrm{ in }\R^n \\
Z_p(0)=|\nabla Z_p(0)|=0 & \\
\lim_{|\x|\to \infty}\frac{Z_p(\x)}{|\x|^3}=0 & \\
\Pi(Z_p,1)=0 & 
\end{array}
\end{equation}
where 
$$
p(\x)=\lim_{r_j\to 0}\frac{\Pi(u,r_j,\x^0)}{\|\Pi(u,r_j,\x^0)\|_{L^2(B_1)}}.
$$

In order to show regularity for the free boundary near a singular point we 
would have to control the limit 
$$
\lim_{r\to 0}\frac{\Pi(u,r_j,\x^0)}{\|\Pi(u,r_j,\x^0)\|_{L^2(B_1)}}.
$$
If one can show that the limit is unique then it follows that the
blow-up 
$$
\lim_{r\to 0} \left(u(r\x+\x^0)/r^2-\Pi(u,r,\x^0)\right)=Z_p
$$
is unique.

The following result, Corollary 7.3 in \cite{ASW3}, gives a 
quantitative measure on how the function $Z_{\Pi(u,r,0)}$ controls
the difference between $\Pi(u,r,0)$ and $\Pi(u,r/2,0)$.

\begin{pro}\label{corgest}
Let $u$ solve (\ref{Eq}) in $B_1\subset \R^n$ and
assume that $\sup_{B_1}|u|\le M$,
$u(0)=|\nabla u(0)|=0$, and that for some $\rho\le \rho_0$ and $r\le r_0$,
$$
\sup_{B_1}|\Pi(u,r)|\ge \frac{1}{\rho}.
$$
Then
$$
\sup_{B_1}|\Pi(u,r/2)-\Pi(u,r)-\Pi(Z_{\Pi(u,r)},1/2)|\le C(M,n,\alpha)(\sup_{B_1}|\Pi(u,r)|)^{-\alpha}
$$
for each $\alpha<1/4$.
\end{pro}

In order to estimate $\sup_{B_1(0)}|\Pi(u,r,0)-\Pi(u,r/2,0)|$ we thus need
to be able to calculate $\Pi(Z_{\Pi(u,r,0)},1/2,0)$. We will do this with 
the help of the following theorem from \cite{KarpMargulis}

\begin{thm}\label{KarpM}
Let $\sigma\in L^\infty(\mathbb{R}^n)$ be homogeneous of
zeroth order, that is $\sigma(\x)=\sigma(r\x)$ for all $r>0$. Assume that
$\sigma$ has the Fourier series expansion
$$
\sigma(\x)=\sum_{i=0}^\infty a_i \sigma_i,
$$
on the unit sphere, where $\sigma_i$ is a homogeneous harmonic polynomial of order
$i$.

Moreover assume that
$\Delta Z= \sigma$ and that
$Z(0)=|\nabla Z(0)|=\lim_{\x\to \infty}Z(\x)/|\x|^3=0$. Then
\begin{align*}
&Z(\x)=q(\x)\ln|\x|+|\x|^2\phi(\x),
\intertext{where}
&q=\frac{a_2}{n+2}\sigma_2
\intertext{and}
&\phi(\x)=
\sum_{i\ne 2} \frac{a_i}{(n+i)(i-2)}\sigma_i\big(\frac{\x}{|\x|}\big).
\end{align*}
\end{thm}

Our strategy in the rest of the paper will be to use Theorem \ref{KarpM}
to calculate 
\begin{equation}\label{SomethingWithZ}
\Pi(Z_{\Pi(u,r,0)},1/2,0)=-\frac{\ln(2)a_2}{n+2}\sigma_2(\x)
\end{equation}
where $\sigma_2$ is the second order term in the Fourier series expansion
$$
-\chi_{\{\Pi(u,r,0)>0\}}=\sum_{i=0}^\infty a_i \sigma_i(\x) \textrm{ on }\partial B_1(0).
$$
Using the expression (\ref{SomethingWithZ}) in Proposition 
\ref{corgest} will give us enough information to 
deduce  that the blow-up of $u$ is unique at all points $\x^0\in S_{n-2}(u)$.

\section{Estimates of the Projections.}

In order to estimate $\Pi(Z_p,1/2)$ we need to calculate $a_2\sigma_2$ from
Theorem \ref{KarpM}. That involves calculating the second order Fourier
coefficients for $-\chi_{\{p_\delta>0\}}$ on the unit sphere. To that end we 
choose $n x_i^2-|\x|^2$ for $i=1,...,n$ and $x_ix_j$ for $i\ne j$ as
a basis for the second order harmonic polynomials. 

We may choose coordinates so that
\begin{equation}\label{pdeltadef}
\frac{\Pi(u,r,0)}{\sup_{B_1}|\Pi(u,r,0)|}=p_\delta(x)=\delta_1 x_1^2+ \delta_2 x_2^2+...+\delta_{n-2} x_{n-2}^2+(1-\tilde{\delta})x_{n-1}^2-x_n^2,
\end{equation}
where $\delta=(\delta_1,\delta_2,...,\delta_{n-2})$ and 
$\tilde{\delta}=\sum_{i=1}^{n-2}\delta_i$. We also define the 
polynomial $p_\delta$, for a given vector $\delta\in \R^{n-2}$ in equation
(\ref{pdeltadef}). We will assume, for definiteness 
that $\tilde{\delta}\ge 0$
(this is implicit in the definition of $p_\delta$ in equation 
(\ref{pdeltadef})). If $\tilde{\delta}<0$ then all the following arguments
follows through with minor and trivial changes.

It follows from symmetry (i.e. $-\chi_{\{p_\delta>0\}}$ is even and the 
$x_i x_j$'s are odd on the unit sphere) that the Fourier coefficient of 
$x_i x_j$ is zero.

Since we are only interested in points $\x^0\in S_{n-2}(u)$ 
where
$$
\lim_{r_j\to 0}\frac{\Pi(u,r_j,\x^0)}{\sup_{B_1}|\Pi(u,r_j,\x^0)|}=p_0,
$$
for some sequence $r_j\to 0$, we may assume 
that $|\delta|< < 1$. 

We also denote by $B_i(\delta)$ the following integral
\begin{equation}\label{Bi}
B_i(\delta)=-\int_{\partial B_1(0)}\chi_{\{p_\delta>0\}}x_i^2 dA
\end{equation}
and by $B(\delta)$ the following integral
\begin{equation}\label{B}
B(\delta)=-\int_{\partial B_1(0)}\chi_{\{p_\delta>0\}} dA.
\end{equation}
Here $dA$ is the surface element.
It follows that the Fourier coefficient of $nx_i^2-|\x|^2$ of $\chi_{\{p_\delta>0\}}$ is
$$
\frac{1}{\|n x_i^2-|\x|^2\|_{L^2(\partial B_1(0))}}\left(nB_i(\delta)-B(\delta)\right).
$$

Using that $\Pi(Z_{p_\delta},1)=0$ by definition and Theorem \ref{KarpM}
we may deduce that
\begin{equation}\label{Zprep}
\Pi(Z_{p_\delta},1/2)=-K_0\sum_{i=1}^n\left(n^2B_i(\delta)-nB(\delta) \right)x_i^2,
\end{equation}
where 
$$
K_0=\frac{\ln(2)}{(n+2)\|nx_i^2-|\x|^2\|_{L^2(\partial B_1(0))}}.
$$

It is clear that we need to estimate the functions $B_i(\delta)$ and $B(\delta)$
in order to estimate 
$$
\Pi(u,r)-\Pi(u,r/2)=\Pi(Z_{p_\delta},1/2)+O\left(\|\Pi(u,r)\|_{L^\infty(B_1(0))}^{-\alpha}\right),
$$
where the above equality is a direct consequence of Proposition \ref{corgest}.

Before we can estimate the integrals in (\ref{Bi}) and (\ref{B}) we need 
to introduce some notation for integration on the unit sphere.
We parametrise the unit sphere in $\R^2$ according to
$$
\partial B_1(0)=\left\{\bar{\xi}_1(\phi);\; \phi\in (0,2\pi) \right\},
$$
where $\bar{\xi}_1(\phi)=(\cos(\phi), \sin(\phi))$. Inductively we define, for
$k\ge 2$, the polar coordinates
$$
\bar{\xi}_k(\phi,\psi_1,\psi_2,..., \psi_{k-1})=(\sin(\phi_{k-1})\bar{\xi}_{k-1}(\phi,\psi_1,...,\psi_{k-2}), \cos(\psi_{k-1})).
$$

The unit sphere in $\R^k$ is then defined by
$$
\partial B_1(0) =\left\{\bar{\xi}_{k-1}(\phi,\psi_1,...,\psi_{k-2});\; \phi\in (0,2\pi),\; \psi_j\in (0,\pi)  \right\},
$$
modulo a set of measure zero.

With this parametrisation an area element on the unit sphere 
becomes
\begin{equation}\label{areael}
dA=\textrm{det}
\left[
\begin{array}{lllll}
\frac{\partial \xi_{k-1}}{\partial \phi} & \frac{\partial \xi_{k-1}}{\partial \psi_1} & \frac{\partial \xi_{k-1}}{\partial \psi_2} & \dots & \frac{\partial \xi_{k-1}}{\partial \psi_{k-2}}
\end{array}
\right]
d\phi d\psi_1...d\psi_{k-2}, 
\end{equation}
where $\xi_{k-1}$ is considered to be a column vector. Somewhat more explicitly
the $k\times (k-1)-$matrix in (\ref{areael}) is
{\small
\begin{equation}\label{matrix}
\left[
\begin{array}{lllll}
-\sin(\psi)P_{j=1}^{k-2} &  \cos(\phi_1)\cos(\phi)P_{j=2}^{k-2} & \dots  & \cos(\psi_{k-3})\cos(\phi)P_{j=1, j\ne k-3}^{k-2} & \cos(\psi_{k-2})\cos(\phi)P_{j=1}^{k-3} \\
\cos(\phi)P_{j=1}^{k-2} & \sin(\phi_1)\cos(\phi)P_{j=2}^{k-2} & \dots &  \cos(\psi_{k-3})\sin(\phi)P_{j=1, j\ne k-3}^{k-2} & \cos(\psi_{k-2})\sin(\phi)P_{j=1}^{k-3} \\
0 & -\sin(\psi_1)P_{j=2}^{k-2} & \dots & \cos(\psi_{k-3})P_{j=1, j\ne k-3}^{k-2} & \cos(\psi_{k-2})P_{j=1}^{k-3} \\
0 &  0 & \ddots & \cos(\psi_{k-3})P_{j=2, j\ne k-3}^{k-2} & \cos(\psi_{k-2})P_{j=2}^{k-3} \\
 \vdots  & \vdots & \ddots  & \vdots  & \vdots \\
0  & 0 & \dots & -\sin(\psi_{k-3})\sin(\psi_{k-2}) & \cos(\psi_{k-2})\cos(\psi_{k-3})  \\
0  & 0 & \dots & 0 & -\sin(\psi_{k-2}) 
\end{array}
\right]
\end{equation}
}
where we have used the notation
$$
P_{j=1}^{k-2}=\Pi_{j=1}^{k-1}\sin(\psi_j).
$$

We will denote the matrix in (\ref{matrix}) by $M$. By the anti-commutativity
of the rows in the determinant function we have the identity
$$
det(M)=
\sin^{k-2}(\psi_{k-3})\sin^{k-2}(\psi_{k-2})\times
$$
\begin{equation}\label{split}
\times\left|
\begin{array}{llllll}
 & & & & 0 & 0 \\
 & & N  & & 0 & 0 \\
 & & & & 0 & 0 \\
0  & 0 & 0  & 0 & -\sin(\psi_{k-3})\sin(\psi_{k-2})  & \cos(\psi_{k-2})\cos(\psi_{k-3}) \\
0 & 0 & 0 & 0 & 0 & -\sin(\psi_{k-2}) \\
\end{array}
\right|
\end{equation}
$$
=\sin^{k-3}(\psi_{k-3})\sin^{k-2}(\psi_{k-2}) det(N)
$$
where $N(\phi,\psi_1,..,\psi_{k-4})$ is the $(k-2)\times (k-3)-$matrix 
satisfying $\sin(\psi_{k-3})\sin(\psi_{k-2})n_{ij}=m_{ij}$ for $1\le i \le k-2$ and $1\le j\le k-3$. Notice that $N$ is independent of $\psi_{k-3}$ and
$\psi_{k-2}$.

In order to estimate $B_i$ we will use the identity in (\ref{split})
to write, with $k=n$,
$$
B_i(\delta)=-\int_{\partial B_1(0)}\chi_{\{p_\delta>0\}}x_i^2 dA_{\partial B_1(0)}=
$$
\begin{equation}\label{IntegralSplit}
-\int_0^{2\pi}\int_0^\pi \dots \int_0^\pi \chi_{\{p_\delta>0\}}x_i^2|det(M)|d\psi_{k-2}d\psi_{k-3}...d\phi=
\end{equation}
$$
-\int_0^{2\pi}\int_0^\pi \dots \int_0^\pi\left[ \int_0^\pi\int_0^\pi \chi_{\{p_\delta>0\}}x_i^2\left|\sin^{n-1}(\psi_{n-3})\sin^n(\psi_{n-2})\right|d\psi_{n-2}d\psi_{n-3}\right]|det(N)| d\psi_{n-4}...d\phi.
$$

We will need some further simplifications
{\small
$$
B_i(\delta)=-\int_{\partial B_1(0)}\chi_{\{p_\delta>0\}}x_i^2 dA_{\partial B_1(0)}
$$
\begin{equation}\label{kex}
=-\int_0^{2\pi}\int_0^\pi \dots \int_0^\pi \chi_{\{p_\delta>0\}}x_i^2|det(M)|d\psi_{k-1}d\psi_{k-2}...d\phi
\end{equation}
$$
=-2^n\int_{(0,\pi/2)^{n-2}}\left[ \int_0^{\pi/2}\int_0^{\pi/2} \chi_{\{p_\delta>0\}}x_i^2S^{n-1,n}(\psi_{n-2},\psi_{n-1})d\psi_{n-1}d\psi_{n-2}\right]|det(N)| d\psi_{n-3}...d\phi
$$
$$
=-2^n\int_{(0,\pi/2)^{n-2}}\left[ \int_{A(\mu)} \chi_{\{p_\delta>0\}}x_i^2S^{n-1,n}(\psi_{n-2},\psi_{n-1})d\psi_{n-1}d\psi_{n-2}\right]|det(N)| d\psi_{n-3}...d\phi
$$
$$
- 2^n\int_{(0,\pi/2)^{n-2}}\left[ \int_{(0,\pi/2)^2\setminus A(\mu)} \chi_{\{p_\delta>0\}}x_i^2S(\psi_{n-2},\psi_{n-1})d\psi_{n-1}d\psi_{n-2}\right]|det(N)| d\psi_{n-3}...d\phi
$$
$$
=I_{1,i}(\delta,\mu)+I_{2,i}(\delta,\mu),
$$
}
where 
$$
S^{n-1,n}(\psi_{n-2},\psi_{n-1})=\left|\sin^{n-1}(\psi_{n-2})\sin^n(\psi_{n-1})\right|,
$$
and $A(\mu)=F^{-1}((0,\mu)^2)$ where $F$ is the stereographic projection
\begin{equation}\label{Fdef}
F(\psi_{n-3},\psi_{n-2})=\left( \frac{\cos(\psi_{n-3})}{\sin(\psi_{n-3})},
\frac{\cos(\psi_{n-2})}{\sin(\psi_{n-2})\sin(\psi_{n-3})}\right).
\end{equation}
If $\mu$ is small then $A(\mu)\approx (\pi/2-\mu,\pi/2)^2$, the exact form of
$A(\mu)$ is unimportant as long as $A(\mu)$ contains a small neighbourhood
of the point $(\pi/2,\pi/2)$. We choose the particular form of $A(\mu)$
in order to simplify some calculations further on (see equation (\ref{somechangeofvar})).

We will estimate $I_{1,i}(\delta,\mu)$ and $I_{2,i}(\delta,\mu)$ separately for 
$|\delta|$ small. Fix a $\mu>0$ such that $|\delta|< < \mu < < 1$. The value of 
$\mu$ is not very important and can be chosen universal, depending only on $n$
in particular $\mu< c_L$ in (\ref{CLdef}).

To estimate $I_{2,i}(\delta,\mu)$ we notice that 
$$
\nabla p_{\delta}=2(\delta_1 x_1, \delta_2 x_2, ..., \delta_{n-2} x_{n-2}^2, (1-\tilde{\delta})x_{n-1}, -x_n ).
$$
By our choice of polar coordinates we have that when 
$\psi_{n-1}\in (0, \pi/2-\mu)$ then 
$$
x_{n}=\cos(\psi_{n-1})\ge c \mu.
$$
This means that the gradient of $p_\delta$ is bounded from below by
a constant times $\mu$ on its zero level set. It is therefore
very easy to estimate $I_{2,i}(\delta,\mu)$ by means of the co-area formula. 

By the co-area formula it follows that for $t\in (0,1)$ and with the notation
$q_\delta=\sum_{j=1}^{n-2}\delta_jx_j^2$
$$
|\frac{d}{dt}I_{2,i}(\delta,t\delta)|
$$
$$
=\left|
\int_{\{(x_{n-2}^2-x_n^2)/(q_\delta-\tilde{\delta}x_{n-1}^2)=t\}}\frac{1}{\left| 
\nabla \frac{x_{n-2}^2-x_n^2}{q_\delta-\tilde{\delta}x_{n-1}^2}
\right|}dA_{\{(x_{n-2}^2-x_n^2)/(q_\delta-\tilde{\delta}x_{n-1}^2)=t\}}
\right| \le C\frac{|\delta|}{\mu}.
$$
In particular
\begin{equation}\label{star1st}
|I_{2,i}(\delta,\mu)-I_{2,i}(\delta,0)|\le C\frac{|\delta|}{\mu}.
\end{equation}

We need to work a little harder in order to estimate $I_1(\delta,\mu)$. We
begin to prove a simple lemma that will allow us to do some
integrations explicitly module $O(|\delta|)$-terms.

\begin{lem}\label{planeproj}
Let $\phi^0,\psi_1^0,\psi_2^0,...,\psi_{n-4}^0$
be fixed. Furthermore we let $\mu>0$ be a small constant and $1\le i\le n$. 
We use polar coordinates $x_i(\phi,\psi_1,\psi_2,..., \psi_{n-2})$.

We also assume that 
\begin{equation}\label{staron5th}
\sum_{j=1}^{n-2}\delta_j x_j(\phi^0,\psi^0_1,\psi^0_2,...,\psi_{n-4}^0 ,\pi/2,\pi/2)^2\ge 0.
\end{equation}

Then there exist a constant $c>0$ such that
$$
\left(1-c\mu \right)\int_{A(\mu)}x_i(\phi^0,\psi_1^0,...,\psi_{n-2})^2
$$
$$
\times \left(\chi_{\{ p_\delta>0\}}(\phi^0,\psi_1^0,...,\psi_{n-2})-
\chi_{\{ p_0>0\}}(\phi^0,\psi_1^0,...,\psi_{n-2})\right)
S^{n-1,n} d\psi_{n-2}d\psi_{n-1}
$$
$$
\le 
\int_0^{\mu}\int_0^\mu\tilde{x}_i^2\left(\chi_{\{p_\delta(\tilde{x})>0\}}(\phi^0,\psi_1^0,...,\psi_{n-3},\psi_{n-2})-
\chi_{\{ p_0>0\}}(\phi^0,\psi_1^0,...,\psi_{n-2})\right) d\tilde{x}_{n-1}d\tilde{x}_{n}
$$
$$
\le \left(1+c\mu \right)\int_{A(\mu)}x_i(\phi^0,\psi_1^0,...,\psi_{n-3},\psi_{n-2})^2
$$
$$
\times \left(\chi_{\{p_\delta>0\}}(\phi^0,\psi_1^0,...,\psi_{n-3},\psi_{n-2})-
\chi_{\{ p_0>0\}}(\phi^0,\psi_1^0,...,\psi_{n-2})\right)
S^{n-1,n}d\psi_{n-3}d\psi_{n-2},
$$
where
$$
S^{i,j} =\left|\sin^i(\psi_{n-3})\sin^j(\psi_{n-2})\right|,
$$
$$
\tilde{x}_i(\phi,\psi_1,\psi_2,...,\psi_{n-2})=\frac{x_i}{\sqrt{\sum_{j=1}^{n-2} x_j^2}}.
$$
and the set $A$ is the stereographic projection of the two dimensional
spherical area 
$$
\{x(\phi^0,\psi_1^0,...,\psi_{n-3},\psi_{n-2});\; (\phi_{n-3},\phi_{n-2})\in (\pi/-\mu,\pi/2)^2\}
$$
under the projection $x\to \tilde{x}$.
\end{lem}

\noindent
{\bf Remark:} {\sl Assumption (\ref{staron5th}) is non-essential
and only made for definiteness and the result still holds if
$$
\sum_{j=1}^{n-2}\delta_j x_j(\phi^0,\psi^0_1,\psi^0_2,...,\psi_{n-4}^0 ,\pi/2,\pi/2)^2 < 0.
$$
}

{\sl Proof:} It is trivial that $1-c\mu\le \sin(\psi_{n-3})\le 1$ and 
that $1-c\mu\le \sin(\psi_{n-3})\le 1$. Therefore 
\begin{equation}\label{Sijest}
1-c_{i,j}\mu\le S^{i,j}\le 1.
\end{equation}

Use the change of variables 
$$
(\psi_{n-3},\psi_{n-2})\to \left( \frac{\cos(\psi_{n-3})}{\sin(\psi_{n-3})},
\frac{\cos(\psi_{n-2})}{\sin(\psi_{n-2})\sin(\psi_{n-3})}\right)=(\tilde{x}_{n-1},\tilde{x}_{n-2})
$$
in
$$
\int_{A(\mu)}x_i(\phi^0,\psi_1^0,...,\psi_{n-3},\psi_{n-2})^2
$$
$$
\times \left(\chi_{\{p_\delta>0\}}(\phi^0,\psi_1^0,...,\psi_{n-3},\psi_{n-2})-
\chi_{\{ p_0>0\}}(\phi^0,\psi_1^0,...,\psi_{n-2})\right)
S^{n-1,n}d\psi_{n-3}d\psi_{n-2}
$$
\begin{equation}\label{somechangeofvar}
=\int_0^\mu\int_0^\mu x_i(\phi^0,\psi_1^0,...,\psi_{n-3},\psi_{n-2})^2\left(\chi_{\{p_\delta(\tilde{x})>0\}}-\chi_{\{ p_0(\tilde{x})>0\}}\right)S^{n-4,n-2}d\tilde{x}_{n-1}d\tilde{x}_n,
\end{equation}
it is in this change of variables that we use the rather awkward definition
of $A(\mu)$ in order to get a nice area of integration to the right.

Since $\sqrt{\sum_{j=1}^{n-2} x_j^2}=\sin(\psi_{n-3})\sin(\psi_{n-2})$
we may estimate 
\begin{equation}\label{tildexiest}
(1-c\mu)\tilde{x}_i\le x_i\le \tilde{x}_i
\end{equation}

Notice that since 
$$
\sum_{j=1}^{n-2}\delta_j x_j(\phi^0,\psi^0_1,\psi^0_2,...,\psi_{n-4}^0 ,\pi/2,\pi/2)^2\ge 0.
$$ 
the integrand is non-negative so we may use (\ref{Sijest}) and
(\ref{tildexiest}) in (\ref{somechangeofvar})
to deduce the desired estimates.
\qed

\begin{lem}\label{deltaladeltaLemma} Let $|\delta|< < \mu < < 1$. 
Also denote
$$
q_\delta=\sum_{j=1}^{n-2}\delta_j x_j^2
$$
and $\phi^0,\psi_1^0,..., \psi^0_{n-4}$ fixed constants.
Then, for $i=1,...,n-2$,
$$
\int_{(\pi/2-\mu,\pi/2)^2}x_i^2\left(\chi_{\{p_\delta>0\}}(\phi,\psi_1,..., \psi_{n-2})-\chi_{\{p_0>0\}}\right)S^{n-1,n}d\psi_{n-3}d\psi_{n-2}=
$$
$$
-\frac{1+O(\mu)}{4}\frac{q_\delta(\phi^0,\psi_1^0,..., \pi/2,\pi/2)}{1-\tilde{\delta}}
\left| \ln(|q_\delta(\phi^0,\psi_1^0,..., \pi/2,\pi/2)|)\right|+O(|\delta|/\mu+
\mu |\delta||\ln(|\delta|)|)
$$
and for $i=n-1,n$ we have
$$
\int_{(\pi/2-\mu,\pi/2)^2}x_i^2\left(\chi_{\{p_\delta>0\}}(\phi,\psi_1,..., \psi_{n-2})-\chi_{\{p_0>0\}}\right)S^{n-1,n}d\psi_{n-3}d\psi_{n-2}
$$
$$
=O(|\delta|)
$$
\end{lem}
{\sl Proof:} By Lemma \ref{planeproj} it is enough to
prove the estimate for
\begin{equation}\label{someintegral}
\int_0^\mu \int_0^\mu x_i^2\left(\chi_{\{p_\delta(x)\}}-\chi_{p_0(x)}\right) dx_{n-1}dx_{n},
\end{equation}
where $\sum_{j=1}^{n-2}x_j^2=1$.

To simplify notation we will write
$$
\kappa=q_\delta(x).
$$
And we will assume that $\kappa>0$, if $\kappa=0$ then the argument is simple 
and the case $\kappa<0$ is treated analogously.

Notice that 
$$
\chi_{\{p_\delta(\tilde{x})>0\}}=
\left\{
\begin{array}{ll}
1 & \textrm{ if } 0< \tilde{x}_n <\sqrt{\kappa+(1+\tilde{\delta})\tilde{x}_{n-1}^2} \\
0 & \textrm{else}.
\end{array}
\right.
$$
For $i=1,...,n-2$ we may write (\ref{someintegral}) as
$$
\int_0^\mu\left(\sqrt{\kappa+(1-\tilde{\delta})x_{n-1}^2)}-\sqrt{x_{n-1}^2)}\right)\tilde{x}_i^2dx_{n-1}
$$
$$
=\frac{1}{4}\frac{\kappa}{1-\tilde{\delta}}\ln(\kappa)\tilde{x}_i^2+\frac{\mu}{2(1+\mu^2)} +O(|\delta|/\mu+
\mu |\delta||\ln(|\delta|)|),
$$
where we have used the identity
$$
\int \sqrt{1+x^2}dx=\frac{1}{2}x\sqrt{1+x^2}+\frac{1}{2}\ln\left(x+\sqrt{1+x^2} \right)
$$
to evaluate the integral. 

For $i=n-1$ we can calculate
$$
\int_0^\mu\left(\sqrt{\kappa+(1+\tilde{\delta})x_{n-1}^2}-\sqrt{x_{n-1}^2}\right)\tilde{x}_{n-1}^2dx_{n-1}=O(\mu^2\kappa)
$$.

Finally, for $i=n$ we get
$$
\int_0^\mu \int_0^\mu x_n^2\left(\chi_{\{p_\delta(x)\}}-\chi_{p_0(x)}\right) dx_{n-1}dx_{n} 
$$
$$
=\int_0^{\mu}\left[ \int_0^{\sqrt{\kappa+(1-\tilde{\delta}) x_{n-1}^2}}x_n^2dx_n-\int_0^{x_{n-1}}x_n^2dx_n\right]dx_{n-1}=O(\kappa \mu^2).
$$
\qed

\begin{pro}\label{Cjest}
If $|\delta|$ is small enough and $C_i(\delta)$ is defined according to
$$
C_i(\delta)=B_i(\delta)-B_i(0)
$$ 
then there exists a universal constant $c$ such that
$$
\frac{1}{c}|\delta \ln(|\delta|)|\le \sum_{j=1}^{n-2} |C_i(\delta)|\le c|\delta \ln(|\delta|)|.
$$

Moreover, if $\delta_i>\delta_j$ then $C_i(\delta)<C_j(\delta)$.
\end{pro}
{\sl Proof:} In (\ref{kex}) we showed that we can write
$$
B_i(\delta)-B_i(0)=\left[ I_{2,i}(\delta,\mu)-I_{2,i}(0,\mu)\right]+
\left[ I_1(\delta,\mu)-I_1(0,\mu)\right].
$$
We also showed, (\ref{star1st}), that
$$
\left[ I_{2,i}(\delta,\mu)-I_{2,i}(0,\mu)\right]=O(|\delta|/\mu).
$$
Also in (\ref{kex}) we showed that 
we can write
\begin{equation}\label{hata}
I_1(\delta,\mu)-I_1(0,\mu)
\end{equation}
$$
=\int_{B_1^{n-2}}\left[\int_A\left(\chi_{\{p_\delta>0\}}-\chi_{p_0>0} \right)S^{n-1,n}d\psi_{n-3}d\psi_{n-2}\right]det(N)dA_{\partial B_1^{n-2}}(\phi,...,\psi_{n-4}).
$$
Furthermore we showed, in Lemmas \ref{planeproj} and \ref{deltaladeltaLemma}, that the inner 
integral in (\ref{hata}) satisfies
$$
\int_A\left(\chi_{\{p_\delta>0\}}-\chi_{p_0>0} \right)S^{n-1,n}d\psi_{n-3}d\psi_{n-2}
$$
$$
=(1+O(\mu))\int_A x_i^2\left(\chi_{\{p_\delta(\tilde{x})>0\}}-
\chi_{\{ p_0>0\}}\right) dx_{n-1}dx_{n}
$$
$$
= -\frac{1+O(\mu)}{4}q_\delta(x_1,..,x_{n-2})
\left| |\ln(|q_\delta(x_1,...,x_{n-2}|)|\right|+O(|\delta|/\mu+
\mu |\delta||\ln(|\delta|)|)
$$
for $(x_1,...,x_{n-2})\in \partial B^{n-2}_1$. Disregarding lower order terms
we may conclude that
\begin{equation}\label{someI1est}
I_{1,i}(\delta,\mu)-I_{1,i}(0,\mu)
\end{equation}
$$
= -\frac{1}{4}\int_{\partial B_1^{n-2}}q_\delta
|\ln(q_{\delta})|det(N)dA_{\partial B_1^{n-2}}+O(|\delta|/\mu+\mu |\delta||\ln(|\delta|)|).
$$

Let us denote the integrand $F(q_\delta)$, that is $F(t)=t|\ln(|t|)|$. We may 
estimate 
\begin{equation}\label{dJens}
\left|F(q_\delta)-|\delta||\ln(|\delta|)\bar{q}_\delta| \right|\le \left|\delta \bar{q}_\delta\ln(|\bar{q}_\delta|)\right| 
\end{equation}
where $\bar{q}_\delta=\frac{1}{|\delta|}q_\delta$. Since $\bar{q}_\delta$
is a second order polynomial with coefficients bounded by one it directly
follows that
\begin{equation}\label{lowerqdelta}
\left|\int_{\partial B_1^{n-2}}\left|\delta\bar{q}_\delta\ln(|q_\delta|)\right| det(N)dA_{\partial B_1^{n-2}}\right|=O(|\delta|).
\end{equation}
By (\ref{lowerqdelta}), (\ref{dJens}) and (\ref{someI1est}) we may estimate
\begin{equation}\label{shitfuck}
I_{1,i}(\delta,\mu)-I_{1,i}(0,\mu)
\end{equation}
$$
=-\frac{|\ln(|\delta|)|}{4}
\int_{\partial B_1^{n-2}}\bar{q}_\delta
det(N) x_i^2dA_{\partial B_1^{n-2}}+O(|\delta|/\mu+\mu |\delta||\ln(|\delta|)|).
$$

We define the linear functional $L:\, \R^{n-2}\to \R^{n-2}$ by
$$
L\delta=\left[
\begin{array}{l}
\int_{\partial B_1^{n-2}}\bar{q}_\delta  x_1^2dA_{\partial B_1^{n-2}}\\
\vdots \\
\int_{\partial B_1^{n-2}}\bar{q}_\delta x_{n-2}^2dA_{\partial B_1^{n-2}}
\end{array}
\right].
$$
Writing $L$ in matrix form we get
$$
L=\lambda_1 I+ \lambda_2 J
$$
where $\lambda_1,\lambda_2>0$, $I$ is the identity matrix and
$$
J=\left[
\begin{array}{lllll}
1 & 1 & 1 & \dots & 1 \\
1 & 1 & 1 & \dots & 1 \\
\vdots & \vdots & \vdots & \ddots & \vdots \\
1 & 1 & 1 & \dots & 1 \\
1 & 1 & 1 & \dots & 1 
\end{array}
\right].
$$
It is easy to see that $\nu^i=[1,1,1,...,1]^T$ is an eigenvector
corresponding to the eigenvalue $\lambda_1+(n-2)\lambda_2$ and
that $\nu^j=e_1-e_j$ for $j=2,...,n-2$ are eigenvectors corresponding to 
the eigenvalue $\lambda_1$. In particular $L$ have $(n-2)-$linearly independent
eigenvectors that correspond to strictly positive eigenvalues. We
may conclude that det$(L)>0$. It follows that there exist a universal constant
$c_L>0$ such that 
\begin{equation}\label{CLdef}
|L\delta|>c_L|\delta|.
\end{equation}

To finish the proof we notice that
$$
\sum_{j=1}^{n-2}|C_i(\delta)|=\sum_{j=1}^{n-2}|B_i(\delta)-B_i(0)|=
|\ln(|\delta|)| |L\delta|+ O(|\delta|/\mu+\mu |\delta||\ln(|\delta|)|)
$$
$$
>\frac{1}{c}|\delta||\ln(|\delta|)|+ O(|\delta|/\mu+\mu |\delta||\ln(|\delta|)|).
$$ 
And
$$
\sum_{j=1}^{n-2}|C_i(\delta)|=\sum_{j=1}^{n-2}|B_i(\delta)-B_i(0)|=
|\ln(|\delta|)| |L\delta|+ O(|\delta|/\mu+\mu |\delta||\ln(|\delta|)|)
$$
$$
< c|\delta||\ln(|\delta|)|+ O(|\delta|/\mu+\mu |\delta||\ln(|\delta|)|).
$$ 
The proposition follows for $\mu$ small enough if $|\delta|< < \mu$. 

The final statement follows easily since $\lambda_1>0$.
\qed 

{\bf Remark:} {\sl We will also use the notation $C(\delta)=B(\delta)-B(0)$.
Notice that 
\begin{equation}\label{Cdeltasum}
C(\delta)=\sum_{i-1}^nC_i(\delta)
\end{equation}
since $\sum_{i=1}^n x_i^2=1$
on the unit sphere.}

\section{Proof of The Main Theorem.}

In this section we prove Theorem \ref{Main}. 

By assumption we have
$$
\lim_{r_j\to 0}\frac{u(r_jx)}{\|u(r_j x)\|_{L^2(B_1)}}=\frac{x_{n-1}^2-x_n^2}{ \|x_{n-1}^2-x_n^2\|_{L^2(B_1)}}
$$
for some sequence $r_j\to 0$. Therefore
\begin{equation}\label{2star2nd}
\lim_{r_j\to 0}\frac{\Pi(u,r_j,0)}{\sup_{B_1}|\Pi(u,r_j,)\|_{L^2(B_1)}}=x_{n-1}^2-x_n^2.
\end{equation}
For any $r>0$ we can define 
a $\delta(r)$ according to 
$$
\frac{\Pi(u,r,0)}{\sup_{B_1}|\Pi(u,r,0)\|_{L^2(B_1)}}=p_{\delta(r)}(x).
$$
With this notation (\ref{2star2nd}) implies that (see \ref{pdeltadef})
$$
|\delta(r_j)|\to 0
$$
so we may, by choosing $j$ large enough, assume that $\delta(r_j)$ is as
small as we need.

Also, from (\ref{loggrowth}) and (\ref{c1alphaBMO}) we may deduce that
$$
\sup_{B_1(0)}|\Pi(u,r_j,0)|\ge c|\ln(r_j)|
$$
for $j$ large enough. 

If we denote $\sup_{B_1(0)}|\Pi(u,s,0)|=\tau_s\approx c|\ln(s)|$ for $s$ small 
enough and $\tau_{2^{-j}s}$ is increasing in $j$ (Lemma \ref{superquad}). 
Then Proposition \ref{corgest} implies that
\begin{equation}\label{deltainrover2}
\Pi(u, r_j/2, 0) = \Pi(u,r_j,0)+\Pi(Z_{p_\delta},1/2,0)+O(\tau_{r_j}^{-\alpha}).
\end{equation}

The main step in our uniqueness proof for blow-up limits is
\begin{lem}\label{deltasmall}
Let $u$ be a solution to (\ref{Eq}) and assume that $\frac{\Pi(u,r,0)}{\sup_{B_1(0)}|\Pi(u,r,0)|}=p_{\delta(r)}$ for some $\delta(r)$ satisfying 
$|\delta(r)|< \kappa_0$ for some universal $\kappa_0$.

We also assume that 
\begin{equation}\label{Cisum}
\sum_{i=1}^{n}C_i(\delta(r))<0. 
\end{equation}

Then for each $\gamma<1/8$ there exist a constant $C_\gamma$ such that
if
\begin{equation}\label{deltalarge}
\max\left(\delta_1(r),\delta_2(r),...,\delta_{n-2}(r) \right)>C_\gamma \tau_r^{-\gamma}
\end{equation}
then
\begin{equation}\label{deltajclaim}
\frac{\max\left(\delta_1(r/2),\delta_2(r/2),...,\delta_{n-2}(r/2) \right)}{1-\tilde{\delta}(r/2)}>
\frac{\max\left(\delta_1(r),\delta_2(r),...,\delta_{n-2}(r) \right)}{1-\tilde{\delta}(r)}.
\end{equation}
Moreover, if  $\delta_j<0$ and
$$
\delta_j\le \min\left(\delta_1(r/2),\delta_2(r/2),...,\delta_{n-2}(r/2) \right)
$$
then it follows that
\begin{equation}\label{negdelta}
\frac{\delta_j(r/2)}{1-\tilde{\delta}(r/2)}< \frac{\delta_j(r)}{1-\tilde{\delta}(r)},
\end{equation}
provided that (\ref{deltalarge}) holds. 
\end{lem}
{\bf Remark:} {\sl If $\sum_{i=1}^{n}C_i(\delta(r))>0$ a similar result holds 
and the proof goes through with trivial changes.}

{\sl Proof:} From (\ref{deltainrover2}) and (\ref{Zprep}) we can conclude that
the coefficient of the $x_j^2-$term in $\Pi(u,r/2,0)$ is
\begin{equation}\label{jcoeff}
\tau_r \delta_j(r)+K_0\left( n^2B_j(\delta(r))-nB(\delta)\right)+O(\tau_r^{-2\gamma}).
\end{equation}

Next we make the following claim

{\bf Claim:} {\sl For $j=1,...,n-2$ we have  $n^2B_j(0)-nB(0)=0$.}

{\sl Proof of the claim:}  This is easy to verify since we can calculate 
$Z_{p_0}$, and thus $B_i(0)$ explicitly (cf. \cite[Lemma 4.4]{ASW2}):

Define $v:(0,+\infty)\times [0,+\infty)\to \mathbb{R}$ by
$$
v(x_{n-1},x_{n}):=
$$
$$
-4x_{n-1}x_n\log(x_{n-1}^2+x_n^2)+2(x_{n-1}^2-x_n^2)\left( \frac{\pi}{2}-2\arctan\left(\frac{x_n}{x_{n-1}}\right)
\right)-\pi(x_{n-1}^2+x_{n}^2).
$$
Moreover, let
$$
w(x_{n-1},x_n):=\left\{ \begin{array}{ll}
v(x_{n-1},x_n),& x_{n-1}x_n\ge 0, x_{n-1}\ne 0,\\
-v(-x_{n-1},x_n),& x_{n-1}<0, x_n\ge 0,\\
-v(x_{n-1},x_n),& x_{n-1}>0, x_n\le 0,
\end{array}\right.$$
and define
$$
\tilde{Z}_{x_{n-1} x_n}(x_{n-1},x_n) := \frac{w(x_{n-1},x_n)-\pi(x_{n-1}^2+x_{n}^2)+8x_{n-1}x_n}{8\pi}.
$$

In particular, $\tilde{Z}_{x_{n-1} x_n}(x_{n-1},x_n)$ is a rotation of $Z_{p_0}$.
It is clear that 
$$
\Pi(\tilde{Z}_{x_{n-1} x_n},1/2,0)=\frac{\ln(2)}{\pi}x_{n-1}x_n,
$$
or equivalently
$$
\Pi(Z_{p_0},1/2,0)=\frac{\ln(2)}{2\pi}(x_{n-1}^2-x_n^2).
$$
It follows that $n^2B_j(0)-nB(0)=0$ for $j=1,...,n-2$. This proves the
claim.

By the definition of $C_j(\delta)$ we may thus write, for $j=1,...,n-2$, 
the coefficient of the $x_j^2-$term in $\Pi(u,r/2,0)$  (that is equation 
(\ref{jcoeff}))
$$
\tau_r \delta_j(r)-K_0\left( n^2C_j(\delta(r))-nC(\delta)\right)+O(\tau_r^{-2\gamma}).
$$
Similarly we can express the $x_{n-1}^2$  coefficient of $\Pi(u,r/2,0)$
according to
$$
\tau_r (1-\tilde{\delta}(r))+\frac{\ln(2)}{2\pi}-K_0\left( n^2C_{n-1}(\delta(r))-nC(\delta)\right)+O(\tau_r^{-2\gamma}).
$$

The quotient of the $x_j^2$ and the $x_{n-1}^2$ coefficients of $\Pi(u,r/2,0)$
is thus equal to
$$
\frac{\tau_r \delta_j(r)-K_0\left( n^2C_j(\delta(r))-nC(\delta)\right)+O(\tau_r^{-2\gamma})}{\tau_r (1-\tilde{\delta}(r)+\frac{\ln(2)}{2\pi}-K_0\left( n^2C_{n-1}(\delta(r))-nC(\delta)\right)+O(\tau_r^{-2\gamma})}.
$$

Let us first prove the Lemma under the assumption 
\begin{equation}\label{deltajmax}
\delta_j(r)=\max\left(\delta_1(r),\delta_2(r),...,\delta_{n-2}(r) \right).
\end{equation}
Then the claim of the Lemma is
\begin{equation}\label{whatweneed}
\frac{\tau_r \delta_j(r)-K_0\left( n^2C_j(\delta(r))-nC(\delta)\right)+O(\tau_r^{-2\gamma})}{\tau_r (1-\tilde{\delta}(r))+\frac{\ln(2)}{2\pi}-K_0\left( n^2C_{n-1}(\delta(r))-nC(\delta)\right)+O(\tau_r^{-2\gamma})}> \frac{\delta_j(r)}{1-\tilde{\delta}(r)}.
\end{equation}

The inequality (\ref{whatweneed}) hold if 
\begin{equation}\label{whatweneed2}
-K_0(1-\tilde{\delta}(r))n^2C_j(\delta(r))+K_0 n(1-\tilde{\delta}(r)-\delta_j(r))C(\delta)+O\left(|\delta|\delta_j+\tau_r^{-2\gamma}\right)>0.
\end{equation}
From (\ref{deltajmax}) and Proposition \ref{Cjest}
we have
$$
(n-1)C_j(\delta)\le \sum_{i=1}^{n-2}C_i(\delta) +O(|\delta|)
=\sum_{i=1}^{n}C_i(\delta)=C(\delta)
$$
where we used Lemma \ref{deltaladeltaLemma} in the first equality
and (\ref{Cdeltasum}) in the last equality. Using
this and $\delta_j>0$ in (\ref{whatweneed2}) we can deduce that the Lemma 
holds if
$$
-K_0(1-\tilde{\delta})C_j(\delta)>O\left(|\delta|\delta_j+\tau_r^{-2\gamma}\right),
$$
or equivalently if 
$$
-C_j(\delta)>O\left(\tau_r^{-2\gamma}\right),
$$
where we used that $|C_j(\delta)|\approx |\delta||\ln(|\delta|)|$.

In particular if $|\delta|$ is small and (\ref{deltajmax}) holds then 
(\ref{deltajclaim}) holds if $\delta_j\ge C_\gamma \tau^{-\gamma}$. 
This is exactly what we wanted to prove.

$ $

Next we chose any $\delta_j<0$ in order to prove (\ref{negdelta}).

Then the claim of the Lemma is
\begin{equation}\label{whatweneedneg}
\frac{\tau_r \delta_j(r)-K_0\left( n^2C_j(\delta(r))-nC(\delta)\right)+O(\tau_r^{-2\gamma})}{\tau_r (1-\tilde{\delta}(r))+\frac{\ln(2)}{2\pi}-K_0\left( n^2C_{n-1}(\delta(r))-nC(\delta)\right)+O(\tau_r^{-2\gamma})}< \frac{\delta_j(r)}{1-\tilde{\delta}(r)}.
\end{equation}

The inequality (\ref{whatweneedneg}) hold if 
\begin{equation}\label{whatweneedneg2}
-K_0(1-\tilde{\delta}(r))n^2C_j(\delta(r))+K_0 n(1-\tilde{\delta}(r)-\delta_j(r))C(\delta)+O\left(|\delta|\delta_j+\tau_r^{-2\gamma}\right)<0.
\end{equation}

We either have that 
\begin{equation}\label{somecase1}
C(\delta)<-\tilde{C}\tau_r^{-\gamma} 
\end{equation}
or 
\begin{equation}\label{somecase2}
C_j(\delta)< \tilde{C} \tau_r^{-\gamma} 
\end{equation}
for some universal $\tilde{C}$. This 
since if 
$\delta_k=\max\left(\delta_1(r/2),\delta_2(r/2),...,\delta_{n-2}(r/2) \right)\ge C_\gamma\tau_r^{-\gamma}$ 
then $C_k(\delta)<-c C_\gamma\tau_r^{-\gamma}|\ln(\tau_r)|$ so if 
$C(\delta)\ge -C\tau_r^{-\gamma}$ then at least one of $C_l(\delta)$,
for $l=1,..., n-2$, must satisfy $C_l(\delta)>c C_\gamma\tau_r^{-\gamma}|\ln(\tau_r)|> > C_\gamma \tau_r^{-\gamma}$ since $|\delta|< < 1$. 
By the monotonicity of $C_l(\delta)$ it follows
that $C_j(\delta)> \tilde{C} \tau_r^{-\gamma}$.

In either case (\ref{somecase1}) or (\ref{somecase2}) it follows that
(\ref{whatweneedneg2}) holds true. The Lemma follows. \qed

We may now proceed with our proof of the main Theorem. From Lemma
\ref{deltasmall} and (\ref{Pining}) it follows that
\begin{equation}\label{deltahastobesmall}
|\delta(r)|\le C\tau_r^{-\gamma}.
\end{equation}
If not then we have by Lemma \ref{deltasmall} that 
$$
\max\left(\delta_1(r/2),\delta_2(r/2),...,\delta_{n-2}(r/2) \right)>\max
\left(\delta_1(r),\delta_2(r),...,\delta_{n-2}(r) \right)
$$
if 
$$
\max\left(\delta_1(r),\delta_2(r),...,\delta_{n-2}(r) \right)>0
$$
and
$$
\min\left(\delta_1(r/2),\delta_2(r/2),...,\delta_{n-2}(r/2) \right)<
\min\left(\delta_1(r),\delta_2(r),...,\delta_{n-2}(r) \right)
$$
if
$$
\min\left(\delta_1(r),\delta_2(r),...,\delta_{n-2}(r) \right)<0.
$$
Since $\tau_{r/2^k}> \tau_{r/2^l}$ for $k>l$ we may iterate this
and conclude that if (\ref{deltahastobesmall}) is not true
then
$$
\lim_{k\to \infty}\max\left(\delta_1(r/2^k),\delta_2(r/2^k),...,\delta_{n-2}(r/2^k) \right)\ge
\max\left(\delta_1(r),\delta_2(r),...,\delta_{n-2}(r) \right)
$$
and/or
$$
\lim_{k\to \infty}\min\left(\delta_1(r/2^k),\delta_2(r/2^k),...,\delta_{n-2}(r/2^k) \right)\le
\min\left(\delta_1(r),\delta_2(r),...,\delta_{n-2}(r) \right).
$$
This would contradict (\ref{Pining}).

So (\ref{deltahastobesmall}) has to hold. This implies in particular
that
$$
\left| \frac{\Pi(u,r,0)}{\sup_{B_{1}}|\Pi(u,r,0)|}-
\frac{\Pi(u,r/2,0)}{\sup_{B_{1}}|\Pi(u,r/2,0)|}\right|\le C\frac{\tau_r^{-\gamma}}{\sup_{B_1}|\Pi(u,r,0)|}\le C\tau_{r}^{-1-\gamma}.
$$
We may iterate and conclude that
$$
\left| \frac{\Pi(u,r,0)}{\sup_{B_{1}}|\Pi(u,r,0)|}-
\frac{\Pi(u,r/2^k,0)}{\sup_{B_{1}}|\Pi(u,r/2^k,0)|}\right|\le C\sum_{j=1}^k C\tau_{r/2^k}^{-1-\gamma}\le
$$
\begin{equation}\label{convergent}
\le C\sum_{j=1}^k \left(k\ln(2)+\ln(1/r)\right)^{-1-\gamma}
\end{equation}
since $\tau_{r}>c |\ln(r)|$. Since $\gamma>0$ it follows that 
(\ref{convergent}) is convergent and we may directly conclude that
$$
\lim_{r\to 0}\frac{u(r_jx)}{\|u(r_j x)\|_{L^2(B_1)}}
$$
exists. The first claim (\ref{1stclaim}) of Theorem \ref{Main}
follows.

That $$
S\cap B_{r_0}(0)\cap \left\{ x;\; \sum_{i=1}^{n-2}x_i^2\le \eta(x_{n_1}^2+x_n^2)  \right\}
$$
consists of two $C^1$ manifolds intersection at right angles at the origin
is now standard (see Corollary 9.2 or in \cite{ASW3}).

To prove that 
$$
S_{n-2}\cap B_{r_0}(0)
$$
is contained in a $C^1$ manifold of dimension $(n-2)$ for some small $r_0$
we may proceed as in Theorem 12.2 in \cite{ASW3}. This proves
Theorem \ref{Main}.

\bibliographystyle{plain}
\bibliography{NdimCodus}

\end{document}